\documentclass[11pt]{article}
\usepackage{latexsym, amssymb, amsfonts, amscd}
\usepackage{graphics}
\usepackage[dvips]{epsfig}

\newtheorem{thm}{Theorem}[section]
\newtheorem{prop}[thm]{Proposition}
\newtheorem{lemma}[thm]{Lemma}

\newtheorem{dfn}[thm]{Definition}
 
\newtheorem{rmk}[thm]{Remark}

\newcommand{\reals}{\mathbb R}

\newcommand{\call}{{\cal L}}

\newcommand{\cinf}{C^{\infty}}

\def\qed{\rule{2.3mm}{2.3mm}}

\newcommand{\half}{\frac{1}{2}}

\begin{document}

\title{\bf Poisson fiber bundles and coupling Dirac structures}

\author{
 {\bf  A\"{\i}ssa Wade} \thanks{email: wade@math.psu.edu}\\
{\small Department of Mathematics, The Pennsylvania State University} \\
{\small University Park, PA 16802.} \\
  }
\date{}
\maketitle

\begin{abstract}
We give sufficient conditions for the existence of
 a Dirac structure on the total space of a Poisson fiber bundle
endowed with a compatible connection.
We also show that Cartan  and Cartan-Hannay-Berry
 connections give rise to  coupling Dirac structures. 

\end{abstract}

{\small \it Mathematics Subject Classification (2000)}: 53D17,
 55Rxx, 57Rxx

{\small {\it Key words}:  Poisson structure, fiber bundle, Dirac structure.}

\section{Introduction}

Several constructions of symplectic
forms on the total space of a symplectic fiber bundle
 have appeared in the literature, among others,
  Thurston's construction (see \cite{Th76}, \cite{McD-S}).
We also have the method of coupling forms developed by
 Guillemin, Lerman and Sternberg
(see \cite{S77}, \cite{GS82}, \cite{GLS96}).  
Moreover in \cite{GLSW83}, Gotay, Lashof,
 \'Sniatycki and Weinstein
 gave necessary and sufficient conditions for the existence
 of a pre-symplectic form on  the total space of  a symplectic fiber bundle
  which restricts to the symplectic structure along its fibers. 
Symplectic fiber bundles have been extensively studied
 in recent years. They have many applications in gauge theories.

\medskip
 On the other hand,
 Poisson fiber bundles which are generalizations of symplectic fiber bundles,
were considered by Marsden, Montgomery and Ratiu 
 in connection with the study of  moving systems (see \cite{MMR90}). 
Based on Cartan's theory of classical space-times,
 they introduced the notion of a {\em Cartan-Hannay-Berry connection},
 which is an important tool for the study of moving systems such
as the ball in a rotating hoop.
 Various  examples of systems having the Cartan connection as 
 underlying geometric structure  can be found in \cite{MMR90}.
It turns out that  Cartan and  Cartan-Hannay-Berry
connections give rise to coupling Dirac structures
in the sense of Vaisman  \cite{Va05} (see Section \ref{Cartan} below).
This suggests that a natural framework for the study of 
certain moving systems is the setting of coupling Dirac structures. This also 
 motivates the study of the problem of finding  conditions
under which the Poisson structure along the fibers of a
 Poisson fiber bundle endowed with a Poisson-Ehresmann
connection can be extended to a (non-vertical) Dirac structure.
 Our aim is  to investigate that problem.
 Our main results are Theorems \ref{thm1}  and \ref{thm2}.
 In \cite{DaW05},  we give another construction of  a coupling structure
 on the total space of a Poisson fiber bundle 
 extending the Sternberg-Weinstein phase space of particles in
 a Yang-Mills field to the  setting of coupling Dirac structures.

Here is an outline of the paper. 
Section 2  provides the tools that will be used to prove the main 
 results. In Section 3, we establish
 Theorems \ref{thm1} and \ref{thm2}. In Section 4, we show that
 Cartan and Cartan-Hannay-Berry connections  induce 
 coupling Dirac structures.

\section{Basic definitions and results}
All manifolds are assumed to be  paracompact,
 Hausdorff, smooth and connected.
We also assume that all maps between manifolds are smooth.

\subsection{Poisson fiber bundles}
Let $(F, {\cal V}_F)$ be a  finite-dimensional Poisson manifold.
A {\em Poisson fiber bundle}  is a fiber bundle 
$F \to E \stackrel{\pi}\to B$ whose
 structure group preserves the Poisson structure on $F$.
In other words, there is an open cover $(U_i)$ of $E$ and
 diffeomorphisms $\phi_i: \pi^{-1}(U_i) \to U_i \times F$  satisfying
 the properties:

\medskip 
\noindent {\bf 1.} the following diagram commutes

\begin{displaymath}
\begin{array}{ccccc}
\pi^{-1}(U_i) &&\stackrel{\phi_i}  {\longrightarrow}&& U_i \times F\\
& \pi  \searrow&&\swarrow {\rm pr}\\
&& U_i \\
\end{array}
\end{displaymath}

\noindent {\bf 2.} If $b \in U_i \cap U_j$
 then the transition map $\phi_{ij}(b)=\phi_i(b)\circ \phi_j(b)^{-1}$
 is a Poisson isomorphism of $(F, {\cal V}_F)$. 
\medskip
 
Notice that the Poisson tensor on each fiber
 $E_b$  is given by ${\cal V}_b=(\phi_i(b))^{-1}_* {\cal V}_F$.
It is independent of the local trivialization
 map $\phi_i$. 
Consider the vertical sub-bundle 
 $${\rm Vert}=Ker(T\pi) \subset TE.$$
\noindent  There is a vertical  Poisson bivector
field  ${\cal V} \in {\mathfrak X}^2_{\rm Vert}(E)
$ which coincides with the Poisson structure along the fibers,
i.e. $(i_b)_* {\cal V}_b = {\cal V},$
 where $i_b: E_b \to E$ is the injection map.

\subsection{Integrable geometric data}
Let $ E \stackrel{\pi} \to B$  be a  fiber bundle.
An  {\em Ehresmann connection} on $E$ is a smooth sub-bundle
 ${\rm Hor} \subset TE$ such that $TE= {\rm Hor} \oplus {\rm Vert}.$ 
 This is alternatively defined by a bundle projection 
 map $\Gamma : TE \to {\rm Vert}$, i.e. $\Gamma^2_e=\Gamma_e$
 for every $e \in E$. One has ${\rm Hor}=ker \Gamma$.

\begin{dfn}{\rm \cite{MMR90}
 Let  $\pi: E \rightarrow B$   be a  Poisson fiber bundle
 together with its associated vertical Poisson bivector field
 ${\cal V} \in \mathfrak X_{\rm Vert}^2(E)$.
 An Ehresmann connection $\Gamma$ on $E$ is {\em Poisson}  if ${\cal V}$ is
 preserved by parallel transport. i.e.
$${\cal L}_{hor_{\Gamma}(X)} \ {\cal V}=0,\quad {\rm for \ all} \ 
X \in \mathfrak X (B),$$
 \noindent where 
$hor_{\Gamma}(X)$ is the $\Gamma$-horizontal lift of $X$.
}
\end{dfn}

\begin{dfn}{\rm \cite{Vo00}
Let $\pi: E \rightarrow B$ be a fiber bundle.
A triple $({\cal V}, \Gamma,\overline{\mathbb F})$
 formed by a vertical bivector field
 ${\cal V} \in \mathfrak X_{\rm Vert}^2(E)$,
 an Ehresmann connection $\Gamma$, and  a horizontal
 2-form $\overline{\mathbb F} \in \Omega^2(E)$  is  called 
{\em geometric data}. It is said to be 
 {\em integrable} if  the following properties are satisfied:
\medskip 

$\bullet$ ${\cal V}$ is a Poisson tensor, i.e. $[{\cal V}, {\cal V}]=0.$

$\bullet$ $\Gamma$ is a Poisson-Ehresmann connection with respect to  
${\cal V}$.

$\bullet$ The curvature  2-form of  $\Gamma$
 is a Hamiltonian vector field given by:

$${\rm Curv}_{\Gamma}(X,Y)= 
{\cal V}^\sharp \Big(d (\overline{\mathbb F}(hor_{\Gamma}(X)
, hor_{\Gamma}(Y)) \Big), 
 \quad \forall \  X, Y \in {\mathfrak X} (B).$$

$\bullet$ The 2-form $\overline{\mathbb F}$
 is horizontally closed.}

\end{dfn}
\noindent {\bf Remark.}

\noindent  {\bf a)} Define the operator 
$\partial_{\Gamma}: \Omega^k(B) \otimes  \cinf(E) \rightarrow 
\Omega^{k+1}(B) \otimes  \cinf(E)$ by setting
\begin{eqnarray*} \partial_{\Gamma} \alpha(X_0, \dots , X_{k})
&=&\sum_{i=0}^k (-1)^i \call_{hor_{\Gamma}(X)}(\alpha(X_0, \dots , \widehat{X_i} ,
 \dots, X_{k}))\cr  
& & +\sum_{i<j} (-1)^{i+j} 
\alpha([X_i, X_j], X_0, \dots ,\widehat{X_i} , \dots,
 \widehat{X_j} , \dots,  X_{k}).
\end{eqnarray*}

The fact that $\overline{\mathbb F}$ is  horizontally closed
can be alternatively expressed  by the following equation (see \cite{Vo00})
 $$\partial_{\Gamma}\mathbb F  =0,$$ where $\mathbb F$ is the 2-form defined by
\begin{equation}
\label{2-form}
\mathbb F(X,Y)= \overline{\mathbb F}(hor_{\Gamma}(X), hor_{\Gamma}(Y)).
\end{equation}

\noindent {\bf b)} Let $({\cal V}, \Gamma, \overline{\mathbb F})$ be integrable 
 geometric data. In general $\partial_{\Gamma}^2 \ne 0$,
  but  its restriction to the Casimir valued $k$-forms,
denoted by $\partial_{\cal V}: \Omega^k(S) \otimes \
\mathrm{Casim}(E, {\cal V}) \rightarrow
 \Omega^{k+1}(S) \otimes \mathrm{Casim}(E, {\cal V}),$
  satisfies $\partial_{\cal V}^2=0$.

\bigskip
\noindent {\bf c)} Let $({\cal V}, \Gamma, \overline{\mathbb F})$ be integrable 
geometric data on $E \to B$. Every $\Phi \in \Omega^1(B) \otimes  \mathrm{Casim}(E, {\cal V})$
induces new integrable geometric data $({\cal V}, \Gamma, \overline{\mathbb F}')$,
where the new horizontal 2-form is defined by 
$$\overline{\mathbb F}'( hor_{\Gamma}(X), hor_{\Gamma}(Y))=
\overline{\mathbb F} (hor_{\Gamma}(X), hor_{\Gamma}(Y))+
(\partial_{\Gamma} \Phi)(X,Y),$$ 
\noindent for any $X, Y \in \mathfrak X(B)$.  
In this case, we say  these geometric data are equivalent. 
This defines an equivalence relation among the set of all
 integrable geometric data with a fixed vertical
 Poisson structure and a fixed Poisson-Ehresmann connection.

\subsection{Coupling Dirac structures}
\subsubsection{Dirac structures}

Let $N$ be a finite-dimensional manifold. Consider the canonical
symmetric pairing $ \langle \cdot ,\cdot \rangle_+$ on the vector
bundle $TN\oplus T^*N$ defined by
\[
\langle (X_1,\xi_1), (X_2,\xi_2)\rangle_+ = \half \Big ( \xi _1
(X_2)+\xi _2 (X_1) \Big ).
\]
The space of smooth sections of $TN\oplus
T^\ast N$ is endowed with a bilinear operation,
 called the Courant bracket, which is an extension of
  the Lie bracket of vector fields to $TN\oplus T^*N$ 
   defined by
\[
[(X_1,\xi_1), (X_2,\xi_2)]=([X_1 ,X_2], \ {\cal L}_{X_1} \xi _2-
i_{X_2}d\xi _1) , 
\]
\noindent for all $(X_1,\xi_1), (X_2,\xi_2)$ smooth 
 sections of $TN\oplus T^*N$.

\begin{dfn}{\rm \cite{C90}
 An {\em  almost Dirac structure}  on a manifold $N$ is a sub-bundle
 $L$ of $ TN\oplus T^\ast N \to N$ which is maximally isotropic with respect
 to the symmetric pairing $\langle \cdot ,\cdot \rangle_+$.
 If, in addition, the space 
 of sections $L$ is closed under the Courant bracket then $L$ is
  called a {\em Dirac structure} on $M$.}
\end{dfn}
Basic examples of Dirac structures are regular foliations, 
   Poisson  and pre-symplectic structures
 (see \cite{C90}).

\subsubsection{Induced Dirac structures on submanifolds} 
Let $L$ be a Dirac structure on a manifold $N$, $Q$ a submanifold of $N$.
At every point $q \in Q$, one has a maximal isotropic vector space
 $$(L_Q)_q={ L_q \cap (T_qQ  \oplus T_q^*N) \over  L_q \cap
(\{0\} \oplus {\rm Ann}(T_q Q))}.$$

\noindent Using the  map
 $(L_Q)_q \rightarrow T_qQ  \oplus T_q^*Q$ defined by
$(u,v) \mapsto (u, v_{|T_qQ}),$
 one can identify $(L_Q)_q$ with a subspace of
 $T_qQ  \oplus T_q^*Q$. Moreover,
  $L_Q$ defines a smooth sub-bundle of
 $TQ \oplus T^*Q$ if and only if $L_q \cap (T_qQ  \oplus T_q^*N)$
  has constant dimension. The following result was proved in \cite{C90}

\begin{prop}{\rm  \cite{C90}}
\label{induced}
If $L_q \cap (T_qQ  \oplus T_q^*N)$
   has constant dimension  then $L_Q$ defines a Dirac structure on $Q$.
\end{prop}

\begin{dfn}{\rm 
 A Poisson fiber bundle $F \to E \stackrel{\pi}\to B$ is 
 {\em coherent} if there exists a Dirac structure  $L$ on $E$
 whose  restriction to the fibers coincides with the Poisson structure
 along the fibers and such that }
$L \cap ({\rm Vert} \oplus {\rm Ann(Vert)}) =\{0\}$.
\end{dfn}

We have the following result:

\begin{prop}
Every coherent Poisson fiber bundle $\pi: E\to B$ admits a Poisson-Ehresmann
connection.
\end{prop}

\noindent {\em Proof:} Suppose $F \to E \stackrel{\pi}\to B$ is a coherent
 Poisson fiber bundle. Let $L$ be a Dirac structure on $E$ that coincides with the Poisson structure on the fibers
 and such that $L \cap ({\rm Vert} \oplus {\rm Ann(Vert)}) =\{0\}$. Then 
  $L_x \cap ({\rm Vert}_x \oplus T_x^*E)$ has constant dimension $n=dim F$.
In fact, $L_x \cap ({\rm Vert}_x \oplus T_x^*E)$ is isomorphic
 to $T_x^*E_x$ since the restriction of $L$ to $E_x$
 is the graph of the Poisson bivector field ${\cal V}_{\pi(x)}$. Set
$$H_x(L)=\{Y_x \in T_xE \ | \  \exists \beta_x \in {\rm Ann(Vert}_x)
\ {\rm such \ that } \ (Y_x, \beta_x) \in L_x \}.$$
We have 
$$H_x(L)\cong {\rm Ann}\Big(pr_2(L_x \cap ({\rm Vert}_x \oplus T_x^*E))\Big).$$
It follows that $$dimH_x(L)= dim E -dim(E_x).$$
Hence $$T_xE= {\rm Vert}_x \oplus H_x(L),$$
\noindent for all $x \in E$. This shows that the distribution Hor$(L)$
 defined by the subspaces $H_x(L)\subset T_xE$  is normal to 
 the sub-bundle Vert. We will
 prove that Hor$(L)$ is smooth. Fix a  neighborhood $U$ of
a point $x \in E$ and let  $(Z_i, \eta_i)$, $(X_j, \alpha_j)$ 
be local bases on $U$ for $L$ and 
$L\cap ({\rm Vert} \oplus T^*E)$, respectively. 
A vector $Y$ tangent to the distribution  Hor$(L)$ has the form
$Y= \sum_i f_i Z_i$ with $\langle Y,  \alpha_j \rangle=0$,
 for all $j$.  The existence of smooth solutions for such a system of equations
 implies the smoothness of Hor$(L)$. Consequently, there is an Ehresmann
$\Gamma_L$  connection  associated with Hor$(L)$.
The fact that $\Gamma_L$ is Poisson  is an immediate consequence
of the  integrability of $L$, i.e. the sections of $L$
are closed under the Courant bracket.

\hfill \qed

\begin{dfn}{\rm 
Suppose the geometric data 
 $({\cal V}, \Gamma,\overline{\mathbb F})$
  defined on the fiber bundle
  $E \rightarrow B$ is integrable. Set
\begin{equation}
\label{Equation 2}
L= \Big\{ (X, i_X\overline{\mathbb F})+
 ({\cal V}^\sharp{\alpha}, \alpha) \ | \ 
 X \in {\rm Hor}_{\Gamma}, \ \alpha \in 
{\rm Ann(Hor}_{\Gamma}) \Big\}.
\end{equation}
Then $L$ is called a {\em coupling Dirac structure}. 
}\end{dfn}

We refer the reader to  \cite{Va05} for 
 a more general definition of a coupling Dirac structure 
 on a foliated manifold. Coupling Dirac structures
 naturally appeared in \cite{DuW04}
 when we considered the transverse Poisson structure
 at a presymplectic of a Dirac manifold.

\begin{rmk}
\label{Poisson}

{\bf a.}The Dirac structure $L$ defined  in  (2)
 satisfies $L \cap (TE \oplus \{0\}) =\{0\}$ if and only if 
$\overline{\mathbb F}$ is non-degenerate. In other words, 
$L$ is the graph of a Poisson bivector field if and only if 
$\overline{\mathbb F}$ is non-degenerate.
\medskip

{\bf b.} The distribution ${\cal D}$ given by the 
 set of all horizontal vector fields $X$ satisfying 
$i_X\overline{\mathbb F}=0$ is integrable. It defines a foliation
 ${\cal F}$, called the  characteristic foliation or
 the null foliation of $L$. 
 Moreover, $E/{\cal F}$ is a Poisson manifold  when
$L$ is reducible (see \cite{LWX98}).

\hfill \qed
\end{rmk}

 We have the following result:

\begin{thm}
\label{1-1}
 Let $E \to B$ be a  fiber bundle.  
 The integrability of geometric data
 $({\cal V}, \Gamma,\overline{\mathbb F})$ 
is equivalent to the fact
 that the space of smooth sections of the corresponding
 subbundle $L \subset TE \oplus T^*E$  (defined as in Equation (2)) 
is closed under the Courant bracket.
\end{thm}
\noindent{\em Proof:}  
Consider geometric data $({\cal V}, \Gamma,\overline{\mathbb F})$ 
 and  define its corresponding almost Dirac structure as in Equation (2).
Set $$e_{_X}= \Big(hor_{\Gamma}(X), \ i_{hor_{\Gamma}(X)}\overline{\mathbb F}
\Big) \quad \quad {\rm and} \quad  \quad 
\quad e_{\alpha}= ({\cal V}^{\sharp} (\alpha), \alpha),$$
\noindent for all $X \in \mathfrak X(B)$ and for all
 $\alpha \in $Ann(Hor$_{\Gamma})$. Since
$${\rm Curv}_{\Gamma}(X,Y)= hor_{\Gamma}([X,Y])
- [hor_{\Gamma}(X), hor_{\Gamma}(Y)],$$ 
we get 
\begin{eqnarray*}
[e_{_X}, e_{_Y}]&=&\Big(hor_{\Gamma}([X,Y])-{\rm Curv}_{\Gamma}(X,Y),  \ \
 {\cal L}_{hor_{\Gamma}(X)}(i_{hor_{\Gamma}(Y)}\overline{\mathbb F})-
 i_{hor_{\Gamma}(Y)}d(i_{hor_{\Gamma}(X)}\overline{\mathbb F} )\Big)\cr
&=& \Big( hor_{\Gamma}([X,Y])-{\rm Curv}_{\Gamma}(X,Y), \ \
i_{[hor_{\Gamma}(X), hor_{\Gamma}(Y)]}\overline{\mathbb F} -
d(\overline{\mathbb F}(hor_{\Gamma}(X), hor_{\Gamma}(Y)))\Big).
\end{eqnarray*}
There follows
$$
 \langle [e_{_X}, e_{_Y}], e_\alpha \rangle_+ =
\half \Big\langle  {\cal V}^{\sharp} 
 \Big(d(\overline{\mathbb F}(hor_{\Gamma}(X), hor_{\Gamma}(Y)))\Big)
- {\rm Curv}_{\Gamma}(X,Y), \ \  \alpha 
\Big \rangle_+ ,$$
\noindent for any $X, Y \in \mathfrak X(B)$ and for any
 $\alpha \in $ Ann(Hor$_{\Gamma})$. Hence 
\begin{equation}
\label{Eq 3}
 \langle [e_{_X}, e_{_Y}], e_\alpha \rangle_+=0, \ \forall e_{\alpha}
 \quad \iff  {\rm Curv}_{\Gamma}(X,Y)=
 {\cal V}^{\sharp} 
 \Big(d(\overline{\mathbb F}(hor_{\Gamma}(X), hor_{\Gamma}(Y))\Big).
\end{equation}
\noindent  Moreover, we have
\begin{equation}
\label {Eq 4}
\langle [e_{_X}, e_{_Y}], e_{_Z} \rangle =0 \quad
 \iff  \quad  d\overline{\mathbb F}\Big(hor_{\Gamma}(X), hor_{\Gamma}(Y),
hor_{\Gamma}(Z)\Big)=0 \end{equation}
\noindent for all $X, Y, Z \in \mathfrak X(B)$.
We also have 
\begin{eqnarray*}
[e_{\alpha}, e_{\beta}]&=&\Big([{\cal V}^{\sharp} (\alpha),
 {\cal V}^{\sharp}(\beta)], \
 {\cal L}_{{\cal V}^{\sharp} (\alpha)} \beta -i_{{\cal V}^{\sharp}(\beta)}
 d \alpha  \Big) \cr
&=& \Big({\cal V}^{\sharp} ({\cal L}_{{\cal V}^{\sharp} (\alpha)} \beta
 -i_{{\cal V}^{\sharp} (\beta)}d \alpha)+ 
[{\cal V}, {\cal V}](\alpha, \beta, \cdot), \ \
{\cal L}_{{\cal V}^{\sharp} (\alpha)} \beta -i_{{\cal V}^{\sharp} (\beta)}
 d \alpha  \Big)
\end{eqnarray*}

\noindent for all $\alpha , \beta \in$ Ann(Hor$_{\Gamma})$.
Therefore, $[e_{\alpha}, e_{\beta}]$ is a smooth section of $L$
 if and only if $[{\cal V}, {\cal V}](\alpha, \beta, \cdot)=0$
for all vertical 1-forms $\alpha, \beta$. Since the trivector field
  $[{\cal V}, {\cal V}]$ is vertical this is equivalent to say that 
 all brackets $[e_{\alpha}, e_{\beta}]$ are smooth  sections of $L$
 if and only if
 
\begin{equation}
\label{Eq 5} 
[{\cal V}, {\cal V}]=0. \end{equation}

\medskip
Furthermore, we have 
$$[e_{_X}, e_{\alpha}]=\Big([hor_{\Gamma}(X), \ {\cal V}^{\sharp} (\alpha)], \ 
{\cal L}_{hor_{\Gamma}(X)} \alpha- i_{{\cal V}^{\sharp} (\alpha)} 
d(i_{hor_{\Gamma}(X)}\overline{\mathbb F} \Big),$$
for any  $X  \in \mathfrak X(B)$, $ \alpha \in $ Ann(Hor$_{\Gamma})$.
 Using the fact that $$[hor_{\Gamma}(X), {\cal V}^{\sharp} (\alpha)]=
\Big({\cal L}_{hor_{\Gamma}(X)} \ {\cal V}\Big)(\alpha, \cdot) + 
{\cal V}^{\sharp}({\cal L}_{hor_{\Gamma}(X)} \alpha),$$
one gets
\begin{equation}
\label{Eq 6}
\langle [e_{_X}, e_{\alpha}], \ e_{\beta}\rangle =0 , \quad \forall \alpha , 
\beta \in {\rm Ann(Hor}_{\Gamma}) \quad
\iff  \quad {\cal L}_{hor_{\Gamma}(X)} \ {\cal V}=0. \end{equation}
Relations (3)-(6) show that if $L$ is a Dirac structure then
 $({\cal V}, \Gamma,\overline{\mathbb F})$ is integrable. 
 The converse is true because of  (3)-(6) and 
 the fact that
\begin{eqnarray*}
\langle [e_{_X}, e_{\alpha}], \ e_{_Y} \rangle_+&=&
 i_{hor_{\Gamma}(Y)}i_{hor_{\Gamma}(X)} d \alpha
-{\cal V}^{\sharp}\Big( \alpha,
d(\overline{\mathbb F}(hor_{\Gamma}(X), hor_{\Gamma}(Y)\Big)\cr
&=&-\Big \langle [hor_{\Gamma}(X), hor_{\Gamma}(Y)] +
 {\cal V}^{\sharp}\Big(d(\overline{\mathbb F}(hor_{\Gamma}(X),
 hor_{\Gamma}(Y))\Big), \ \ \alpha \Big \rangle_+ \cr
&=& \Big \langle {\rm Curv}_{\Gamma}(X,Y)- {\cal V}^{\sharp} 
 (d(\overline{\mathbb F}(hor_{\Gamma}(X), hor_{\Gamma}(Y)))
, \ \ \alpha
\Big \rangle_+.
\end{eqnarray*}
\noindent This completes the proof of Theorem \ref{1-1}.

\hfill \qed

\section{Dirac extensions of Poisson fiber bundles}
 In this section, we give  constructions
 of Dirac structures on the total space of certain 
Poisson fiber bundles. First, we recall from 
 \cite{We87} the following definition:

\begin{dfn}{\rm 
A {\em classical Yang-Mills-Higgs} setup is a triple 
 $(G, P, F)$ formed by a finite-dimensional Lie group $G$,
a principal $G$-bundle $P$, and a  Hamiltonian Poisson $G$-space $F$.}
\end{dfn}

\subsection{Dirac structures and principal bundles}
We have the following result:

\begin{thm}
\label{thm1}
Let $(G, P, F)$ be a classical Yang-Mills-Higgs setup. 
Then every   connection $\Theta$ on $P$ gives rise to a coupling 
 Dirac structure on the associated bundle $E=P \times_G F$. 
\end{thm}

\noindent {\em Proof:}
Let $\pi_{_{P \times F}} : P \times F \to P \times_{G} F$ 
denote the canonical projection.
 Define the vertical bivector field ${\cal V}$ on 
 $E$ as  follows 
$${\cal V}= (\pi_{_{P \times F}})_* {\cal V}_F.$$
It satisfies $[{\cal V}, {\cal V}]=0$ since $ {\cal V}_F$ is Poisson.
Moreover, every connection $\Theta$ on $P$  induces a connection $\Gamma$ on
$E$. The $\Gamma$-horizontal lift of $X \in \mathfrak X(B)$ is given by
\begin{equation}
\label{Eq 7}
(hor_{\Gamma}(X))([p,m])=T_{(p,f)}\pi_{_{P \times F}} 
( \overline{X}_p, {\bf 0}_f), \quad \forall \ [p,m] \in P \times_{G} F,
\end{equation}
\noindent where  ${\bf 0}_f$ is the zero tangent vector at $f \in F$ and
$\overline{X}_p\in T_pP$ is the $\Theta$-horizontal lift of $X$ at $p.$
Consequently, one gets
$${\cal L}_{hor_{\Gamma}(X)} \ {\cal V}=0,$$
\noindent for all $X \in \mathfrak X(B)$.
Recall that the curvature of $\Theta$ is a vertical $\mathfrak g$-valued
 2-form. Moreover, for all $X, Y \in \mathfrak X(B)$, we have 
\begin{equation}
\label{Eq 8}
({\rm Curv}_{\Gamma}(X,Y)) ([p, f])= T_{(p, f)}\pi_{_{P \times F}}
\Big( {\bf 0}_p , \ (\varrho_F \circ {\rm Curv}_{\Theta}(\overline{X_p}, 
\overline{Y_p}) )_f \Big),
\end{equation}
\noindent where  $\varrho_F: \mathfrak g \to \mathfrak X (F)$
 is the infinitesimal action associated to the
 $G$-action on $F$. Let $J: F \to \mathfrak g^*$ be the momentum
 map associated with the $G$-action on $F$.
 Now, define  the {\em horizontal } 2-form 
$\overline{\mathbb F}$ as follows

\begin{equation}
\label{curv} 
\Big(\overline{\mathbb F}(hor_{\Gamma}(X),\ hor_{\Gamma}(Y))\Big)([p,f])= 
\Big\langle J(f), \ {\rm Curv}_\Theta (\overline{X_p}, \
 \overline{Y_p}) \Big\rangle,
\end{equation}
\noindent  for all   $X, Y \in \mathfrak X(B)$, and for all
$[p,f] \in E$. Using Relations (8)-(9) and 
the fact that the $G$-action  on 
 $F$  is Hamiltonian, one obtains
$${\rm Curv}_{\Gamma}(X,Y)= 
{\cal V}^\sharp \Big(d (\overline{\mathbb F}(hor_{\Gamma}(X)
, hor_{\Gamma}(Y)) \Big), 
 \quad {\rm for \ all} \  X, Y \in {\mathfrak X} (B).$$
To check that $\overline{\mathbb F}$ is horizontally closed,
 it is enough to notice  that if we  set
$$\Phi_{[p,f]}(\widehat{\mathbf {A}_p, \mathbf {B}_f})=\Big\langle J(f) ,
 \ \Theta_p(\mathbf{A}_p) \Big\rangle , $$
 for all  tangent vectors $$(\widehat{\mathbf{A}_p, \mathbf{B}_f}) =
T_{(p, f)}\pi_{_{P \times F}}(\mathbf{A}_p, \mathbf{B}
_f)$$  then  we get
 $$\overline{\mathbb F}\Big(hor_{\Gamma}(X),\ hor_{\Gamma}(Y)\Big)=
 d \Phi\Big(hor_{\Gamma}(X),\ hor_{\Gamma}(Y)\Big),$$ 
\noindent for all $ X, Y \in {\mathfrak X} (B)$.
The fact that $\overline{\mathbb F}$ is horizontally closed 
follows from $d^2\Phi=0$.
We have constructed integrable geometric data $({\cal V}, \Gamma, 
 \mathbb F)$. Finally, we can apply Theorem  \ref{1-1} which gives
  the result we sought.

\hfill \qed

\subsection{Fat bundles}
Let $P \to B$ be a (left) principal $G$-bundle, ${\mathfrak g}^*$ the
 dual of the Lie algebra of $G$, and $S$ a subset
 of ${\mathfrak g}^*$.
 A connection $\Theta$ on $P$ is {\em $S$-fat} \cite{We80} if for every
 $\mu \in S$,  $\mu \circ {\rm Curv}_{\Theta}$ is non-degenerate.

\begin{prop}
\label{fat}
Let $(G, P, F)$ be a classical Yang-Mills-Higgs setup. 
 Then every $J(F)$-fat connection  $\Theta$  on $P$
gives rise to a Poisson structure
on the associated bundle $E=P \times_G F \to B$. 
\end{prop}

The proof of  Proposition \ref{fat}
 is similar to that of Theorem \ref{thm1}. Precisely,
 one can notice that  the 2-form  defined by Equation  (\ref{curv}) 
 is nondegenerate when the given connection $\Theta$
 is $J(F)$-fat connection. So using Remark \ref{Poisson},
  we conclude that the Dirac structure obtained 
 (as in the proof of Theorem \ref{thm1})
  is the graph of a Poisson bivector field on $E$.

\hfill \qed

\subsection {Another construction of a Dirac extension
 of a Poisson fiber bundle}

\begin{thm}\label{thm2}
Let $(F, {\cal V}_F)$ be a compact Poisson manifold whose first Poisson 
cohomology group $H^1_{{\cal V}_F}(F)$ vanishes.
 Let $F \to E \stackrel{\pi}\to B$  be a  Poisson fiber bundle.
Then every Poisson-Ehresmann connection on $E$
 gives rise to an equivalence class of coupling Dirac structures
 on $E$ such that each representative restricts to 
the Poisson structure on the fibers. 
\end{thm}

\noindent {\em Proof:} Consider the structure
 group $G=$Iso$(F, {\cal V}_F)$ which consists of all
 Poisson isomorphisms of $(F, {\cal V}_F)$. The 
{\em Poisson frame bundle}, denoted by $P$, is the  the principal
 $G$-bundle whose fiber over $b$ is the set of all 
Poisson isomorphisms  $\varphi_b : (F, {\cal V}_F) \to (E_b, {\cal V}_b)$. 
We can identify $E$ with $P \times_G F$. The vertical Poisson
 vector field ${\cal V}$ (defining the Poisson  
 fiber bundle structure) can be viewed as the push-forward 
of $ {\cal V}_F$ by the projection map 
$\pi_{_{P \times F}}: P \times F \to P \times_G F$.
 Moreover, every 
Poisson-Ehresmann connection  $\Gamma$ on $E$ induces a connection
 $\Theta$ on $P$. These connections are related  as in Equation (7).
 
\medskip

  Consider the $\reals$-linear map $J$ from the Lie algebra  of $G$
   into $\cinf(F)/\{\rm Casimir \ functions\}$
 such that $J(Z)=g_{_Z}$ is the unique  function on $F$ 
(up to Casimir functions) whose Hamiltonian vector field equals $Z$.
Notice that the Lie algebra  of $G$ coincides 
with  the space of Hamiltonian  vector fields of $(F, {\cal V}_F)$
  because of the hypothesis $H^1_{{\cal V}_F}(F)=\{0\}$.
Using  this map $J$ and the 
connection 1-form $\Theta(p): T_p P  \to \rm{Ham}(F, {\cal V}_F)$,
 we define a class of 1-forms $\Psi$ on $E$ as follows:
$$\Big(\Psi(Y)\Big)(e)=\Big(J \circ (\Theta(p)(Y^1_p)\Big)(f),$$
\noindent for every $e=[p,f] \in E$ and for all 
$Y \in \mathfrak X(E)$ defined by
 $$Y([p,f])= T_{(p, f)}\pi_{_{P \times F}}(Y^1_p, Y^2_f).$$ 
Define the class of horizontal 2-forms $$
\overline{\mathbb F}\Big(hor_{\Gamma}(X),\ hor_{\Gamma}(Y)\Big)
= d \Psi\Big(hor_{\Gamma}(X),\ hor_{\Gamma}(Y)\Big),$$
 which is determined up to elements of the form
$\overline{\partial_{\Gamma} \Phi}$, where
 $\Phi \in \Omega^1(B) \otimes {\rm Casim}(F, {\cal V}_F)$
and $$\overline{\partial_{\Gamma} \Phi}\Big(hor_{\Gamma}(X),\ hor_{\Gamma}(Y)\Big)=
(\partial_{\Gamma} \Phi)\big(X, Y \big),$$ 
\noindent for $X, Y \in \mathfrak X (B)$. By construction,
 each representative element, also denoted by $\overline{\mathbb F}$, is horizontally closed.
 Furthermore, by arguments similar to those
 used in the proof of Theorem \ref{thm1}, one gets 
$${\rm Curv}_{\Gamma}(X,Y)= 
{\cal V}^\sharp \Big(d (\overline{\mathbb F}(hor_{\Gamma}(X)
, hor_{\Gamma}(Y)) \Big),  \quad {\rm for \ all} \  X, Y \in {\mathfrak X} (B).$$
There follows Theorem \ref{thm2}.

\hfill \qed

\section{The Cartan-Hannay-Berry connection} \label{Cartan}
In this section, our goal is to show that the notion of a Cartan-Hannay-Berry connection
 provides specific examples of coupling Dirac structures.
 We will use the following lemma.

\begin{lemma}
\label{lemme 1}
Let $\pi: F \times B \rightarrow B$  be a Poisson fiber bundle together
with  its associated vertical Poisson bivector field ${\cal V}$.
Consider an Ehresmann connection $\Gamma$ on $E$  such that
$$\Gamma(0, X)= {\cal V}^{\sharp}(d \Phi(X)), \quad \forall \ X \in 
\mathfrak X(B),$$
\noindent for some $\Phi \in \Omega^1(B)\otimes \cinf(E)$.
Set $$\mathbb F (X,Y)= d \Phi(X,Y) -  \{\Phi(X), \Phi(Y)\}_{\cal V}, \quad 
\quad \forall \ X, Y \in 
\mathfrak X(B).$$
Then the curvature of $\Gamma$ is given by 
$${\rm Curv}_{\Gamma}(X,Y)=\Big({\cal V}^{\sharp}
 (d (\mathbb F (X,Y)), \ 0 \Big),$$
\noindent for any $X, Y \in \mathfrak X(B)$.
Moreover, the associated horizontal 2-form $\overline{\mathbb F}$
 (defined as in Equation (1)) is horizontally closed.
\end{lemma}

\noindent The proof of this lemma is straightforward.
It is left to the reader.
Now, we recall from \cite{MMR90}
 the definition and properties of a Cartan connection.

\bigskip 
Let $S$ be a Riemannian manifold, 
 $Q$ the configuration space of a given mechanical system,
 and $B$ a finite-dimensional space of embeddings of $Q$ into $S$.
Given a vector field $U \in {\mathfrak X}(B)$ and a point $b \in B$,
 the tangent vector  $U_b \in T_bB$ is a map 
$U_b: Q \to TS$ with $U_b(q) \in T_{b(q)}S$.
There is a canonical vector field
 ${\cal U}_b \in {\mathfrak X}(Q)$ associated with $U_b$.
It is defined as follows: let $U_b^{\perp}(q)$ be the orthogonal
 projection of $U_b(q)$ to $T_{b(q)} b(Q)$ then 
 $${\cal U}_b(q)= (Tb)^{-1}(U_b^{\perp}(q)).$$
The Cartan connection $\gamma_0$ on the trivial fiber bundle
 $Q \times B \to B$  is given by
$$(\gamma_0(V, U))(q,b)=(V_b+{\cal U}_b(q), 0_b).$$
Consider the Poisson fiber bundle $T^*Q \times B \to B$ with
typical fiber $(T^*Q, \omega_{\rm can})$. Denote by
 ${\cal V} \in {\mathfrak X}^2_{Vert}(T^*Q \times B)$
  the vertical Poisson structure  determined by the Poisson structure
 on the fibers.

\begin{dfn}{\rm \cite{MMR90}
 The {\em induced Cartan connection} on $E=T^*Q \times B $
 is the map $\Gamma_{_0} : TE \to {\rm Vert}$ defined by
$$\Gamma_{_0}(W, U)=(W+X_{{\cal P}(U)}, 0),$$
\noindent where ${\cal P} \in \Omega(B)\otimes \cinf(T^*Q \times B)$
 is the 1-form defining the momentum function of $U$. Precisely,
 we have
$$({\cal P}(U))(\alpha_q, b)= \langle \alpha_q , \ {\cal U}_b(q) \rangle
 \quad \forall \ \alpha_q \in T^*_qQ, \ \forall \ b \in B,$$
and $X_{{\cal P}(U)}$ is the Hamiltonian vector field of ${\cal P}(U)$
 relative to ${\cal V}$. Moreover,
$$(hor_{\Gamma_{_0}}U)= (0, U)+(-X_{{\cal P}(U)},0),
 \quad {\rm for \ every} \ U \in \mathfrak X(B).$$
}
\end{dfn}

Let $G$ be a compact  Lie group.
Given a left action of $G$ on $T^*Q$  with equivariant
 momentum map $J: T^*Q \to \mathfrak g^*$, we denote by 
$\langle \cdot , \cdot \rangle_G$ the averaging operation
 (see \cite{MMR90}).

\begin{dfn}{\rm \cite{MMR90}
 The {\em Cartan-Hannay-Berry connection} on 
 $T^*Q \times B$ is the vertical valued 1-form $\Gamma$ defined as follows:
$$\Gamma (W, U)=(W+X_{ \langle {\cal P}(U) \rangle_G}, 0),$$
for any $W \in \mathfrak X(T^*Q)$, $U \in \mathfrak X(B)$.
In other words, the horizontal lift for $U \in \mathfrak X(B)$ is given by
$$(hor_{\Gamma}U)(\alpha_q, b) =
(-X_{\langle {\cal P}(U) \rangle_G}(\alpha_q, b), \ U(b)).$$
}
\end{dfn}

Now we are going to define the integrable data associated with
 the Cartan-Hannay-Berry connection.
We set
$$\mathbb F_0 (U_1, U_2)= U_1 \cdot {\cal P}(U_2)
- U_2 \cdot {\cal P}(U_1) - {\cal P}([U_1, U_2])-
\{{\cal P}(U_1), {\cal P}(U_2)\}_{\cal V},$$
and $$\mathbb F (U_1, U_2)=
\langle \mathbb F_0 (U_1, U_2) \rangle_G.
$$
\noindent for all $U_1, \ U_2 \in \mathfrak X(B)$.
 Lemma \ref{lemme 1} implies that 
 $({\cal V}, \Gamma_0, \mathbb F_0)$ and
$({\cal V}, \Gamma, \mathbb F)$
 are integrable geometric data on $T^*Q \times B$.
 Their associated coupling Dirac structures 
 are defined as in Equation (1). In other words, every 
 Cartan (resp. Cartan-Hannay-Berry) connection gives rise to a coupling
 Dirac structure.
\hfill \qed

\medskip

\noindent{\bf Acknowledgments.}  I would like to thank Ben Davis,
 Jean-Paul Dufour, Charles-Michel Marle, Jean-Pierre Marco,
  and Alan Weinstein for helpful discussions
and remarks. I am also grateful to the MSRI for its kind hospitality
when the manuscript was being prepared.

\end{document}